# ON UTILITY MAXIMIZATION IN DISCRETE-TIME FINANCIAL MARKET MODELS

By Miklós Rásonyi[1] and Lukasz Stettner[2]

*Computer and Automation Institute of the Hungarian Academy of Sciences and Institute of Mathematics of the Polish Academy of Sciences*

We consider a discrete-time financial market model with finite time horizon and give conditions which guarantee the existence of an optimal strategy for the problem of maximizing expected terminal utility. Equivalent martingale measures are constructed using optimal strategies.

**1. Introduction.** In this paper we study the existence of optimal portfolios for maximizing expected utility at the end of a trading period in a financial market. Preferences of the agent in consideration are described by a nondecreasing concave function $U : \mathbb{R} \to \mathbb{R}$, trading dates occur at discrete time instants.

The same problem has been treated in [27] for a general, continuous-time semimartingale model. The article (similarly to its predecessor [20]) formulated a so-called "reasonable asymptotic elasticity" condition on $U$ which is sufficient for the existence of an optimal stratregy (provided that the asset price process admits an equivalent local martingale measure and is locally bounded). This condition appears to be necessary in the general context, as highlighted by the counterexamples given there. The paper made extensive use of functional analysis and followed an approach via the dual problem. We wondered how this could be avoided, at least in a discrete-time setting.

As it will become clear in the arguments below, in discrete-time market models a direct probabilistic approach is possible, based on a simple

Received November 2003; revised June 2004.
[1]Supported by National Research Foundation of Hungary (OTKA) Grants T 047193 and F 049094.
[2]Supported by PBZ KBN 016/P03/99.

idea going back to, for example, [25]. Under weaker "asymptotic elasticity" conditions, we manage to establish the existence of optimal strategies for nonsmooth utility functions and for possibly unbounded price processes, hence, we cover (though only in discrete-time) several cases where previous results do not apply, see Remark 2.6 of [29], as well as Remark 2.9 below. Although our arguments appear to be fairly straightforward, they are not quite evident due to a number of hidden pitfalls related to some delicate measure theoretic issues.

It seems that, despite its practical importance, discrete-time utility maximization has been somewhat neglected lately. Schäl [31] considers utilities defined on the whole real line in a finite probability space setting; Schäl [29, 30] studies the case of utilities $U:\mathbb{R}_+ \to \mathbb{R}$; Kramkov and Schachermayer [20, 21] give exhaustive treatments of the case $U:\mathbb{R}_+ \to \mathbb{R}$ in a general semimartingale model. It is possible to apply our techniques to this kind of utility function, too; related results will be presented elsewhere.

On the history of the problem, consult the papers cited above with the references therein. In Section 2 below we present our main results. Section 3 deals with consequences of the absence of arbitrage property. Section 4 considers a one-step model which is then carried over to several time steps in Section 5. In Section 6 equivalent martingale measures are constructed using optimal strategies, Section 7 presents examples and corollaries of the main results. Finally, the Appendix contains auxiliary material.

**2. Problem formulation.** A usual setting for discrete-time market models is considered: a probability space $(\Omega, \mathcal{F}, P)$; a filtration $(\mathcal{F}_t)_{0 \leq t \leq T}$ and a $d$-dimensional adapted process $(S_t)_{0 \leq t \leq T}$ describing the (discounted) price of $d$ assets which are present in a given economy. We suppose that $\mathcal{F}_0$ contains all $P$-zero sets. The symbol $\langle \cdot, \cdot \rangle$ denotes the usual scalar product in $\mathbb{R}^d$, $|x| := \sqrt{\langle x, x \rangle}$.

In what follows, $\Xi_t$ will denote the set of $\mathcal{F}_t$-measurable $d$-dimensional random variables. Trading strategies are represented by arbitrary $d$-dimensional predictable processes $(\phi_t)_{1 \leq t \leq T}$, where $\phi_t^i$ denotes the investor's holdings in asset $i$ at time $t$; predictability means that $\phi_t \in \Xi_{t-1}$. The family of all predictable trading strategies is denoted by $\Phi$. In continuous-time models "admissibility" requirements are usually imposed on portfolios (e.g., the value process is bounded from below). An important feature of the present approach is that we go beyond the class of (locally) bounded price processes. When dealing with unbounded $S$, it is crucial to allow portfolios which are not necessarily bounded from below.

The value at time $t$ of a portfolio $\phi$ starting from initial capital $c$ is given by

$$V_t^{c,\phi} = c + \sum_{i=1}^{t} \langle \phi_i, \Delta S_i \rangle,$$

where $\Delta S_i := S_i - S_{i-1}$ and $c \in \mathbb{R}$.

Fix a concave nondecreasing function $U : \mathbb{R} \to \mathbb{R}$. The positive (resp. negative) part of a real-valued function $V$ is denoted by $V^+$ (resp. $V^-$). Regular conditional distributions and generalized conditional expectations are used throughout the paper. Dependence of various functions on $\omega \in \Omega$ will often be dropped in the notation. By convention, $U'(x)$ denotes the *left-hand* derivative of $U$ at $x$.

We are dealing with maximizing the expected terminal utility $EU(V_T^{c,\phi})$ from initial endowment $c$. In order to have a well-posed problem, it should be stipulated that the optimal value is finite.

ASSUMPTION 2.1. Suppose that the following random functions are well defined:
$$U_T(x) := U(x), \qquad x \in \mathbb{R};$$
for $0 \leq t < T$ and for all $x \in \mathbb{R}$,
$$U_t(x) := \operatorname*{ess\,sup}_{\xi \in \Xi_t} E(U_{t+1}(x + \langle \xi, \Delta S_{t+1} \rangle) | \mathcal{F}_t) < \infty \qquad \text{a.s.},$$
and for all $x \in \mathbb{R}$,
$$EU_0(x) < \infty. \tag{1}$$

REMARK 2.2. We remark that regular versions of $U_t$ exist by Proposition 4.4 below. The quantity $U_t(x)$ is the highest future expected conditional utility with respect to $\mathcal{F}_t$ for an agent who starts trading at time $t$ with initial endowment $x$. So Assumption 2.1 roughly says that the supremum of future expected utility at time 0 should not attain $\infty$. This is a natural requirement in the context of utility maximization. If $U$ is bounded from above [e.g., the often encountered exponential utility $U(x) = 1 - e^{-x}$ or the shortfall function $U(x) = \max\{x, 0\}$], then Assumption 2.1 trivially holds. For unbounded $U$, Assumption 2.1 seems to be more difficult to verify, we shall check its validity in a fairly broad model class in Proposition 7.1.

We will impose the following *absence of arbitrage* (NA) property (see Section 3 for a discussion):
$$(NA) : \forall \phi (V_T^{0,\phi} \geq 0 \text{ a.s.} \;\Rightarrow\; V_T^{0,\phi} = 0 \text{ a.s.}). \tag{2}$$

One can assert the existence of an optimal strategy under certain conditions on $U$.

ASSUMPTION 2.3. The utility function $U : \mathbb{R} \to \mathbb{R}$ is concave, nondecreasing; $U(0) = 0$, and there exists $\tilde{x} > 0$ and $0 < \gamma < 1$ such that for $x \geq \tilde{x}$ and for any $\lambda \geq 1$,
$$U(\lambda x) \leq \lambda^\gamma U(x). \tag{3}$$

REMARK 2.4. Condition (2.3) appears as a hypothesis in [20] and [27]; it is equivalent to a certain asymptotic elasticity property of $U$:

$$\limsup_{x \to \infty} \frac{U'(x)x}{U(x)} < 1;$$

see Section 6 of [20] for details. This concept encompasses the most frequently encountered behaviors of utility functions at $\infty$ (bounded, logarithmic, power $< 1$). The condition $U(0) = 0$ can evidently be dropped, we stipulate it only for the sake of a simpler presentation.

It is possible to replace (2.3) by a hypothesis on the behavior of $U$ at $-\infty$.

ASSUMPTION 2.5. The function $U$ is concave, nondecreasing; $U(0) = 0$ and there exist $\alpha > 0, \tilde{x} \leq 0$ such that

(4) $$U(\lambda x) \leq \lambda^{1+\alpha} U(x),$$

for all $x \leq \tilde{x}$.

REMARK 2.6. Again, this is equivalent to another property of a "asymptotic elasticity," which appears in [27]:

$$\liminf_{x \to -\infty} \frac{U'(x)x}{U(x)} > 1,$$

see also [10] and [11].

THEOREM 2.7. *Let $U$ satisfy either Assumption 2.3 or Assumption 2.5 and $S$ satisfy (2). Let us suppose that Assumption 2.1 holds true. Then there exists a strategy $\phi^* = \phi^*(c)$ satisfying*

$$u(c) = EU(V_T^{c,\phi^*}) < \infty,$$

*where*

$$u(c) := \sup_{\phi \in \Phi(U,c)} EU(V_T^{c,\phi}),$$

*and $\Phi(U,c)$ is the set of strategies $\phi \in \Phi$ for which the expectation $EU(V_T^{c,\phi})$ is well defined.*

Introduce the random subset $D_t(\omega)$ of $\mathbb{R}^d$: the smallest affine hyperplane containing the support of the (regular) conditional distribution of $\Delta S_t$ with respect to $\mathcal{F}_{t-1}$; this is an $\mathcal{F}_{t-1}$-measurable random set, see the Appendix and Proposition A.1 in particular. We now present a uniqueness result.

THEOREM 2.8. *If the assumptions of the previous theorem are met and $U$ is strictly concave, then there is a unique optimal strategy $\phi^*$ satisfying*

$$\phi_t^* \in D_t \qquad a.s.$$

We present the proofs of Theorems 2.7 and 2.8 in Sections 4 and 5, under Assumption 2.3. At the end of Section 5 we indicate the necessary modifications under Assumption 2.5.

REMARK 2.9. In [27] the existence of an optimal strategy is proved for locally bounded $S$ (which is a local martingale under some equivalent measure) and under certain conditions on $U$ [differentiability, Inada conditions, (2.3) and (2.5)].

In the present discrete-time setting one may assert the existence of an optimal strategy for a larger class of utility functions. Examples 7.4 and 7.5 show that there are applications of our main result for which this extension is crucial.

We also allow $S$ to be possibly unbounded. In this case, the usual duality approach (see, e.g. [17, 27] or [9]) does not work: according to Counterexample 2.1 on page 46 of [9], the dual problem may fail to admit an optimal solution.

One may even consider even random utility functions $U(x,\omega)$; see [19] and [1]. This comes in handy, for example, when we have a fixed random variable $B(\omega)$ (contingent claim) and try to maximize

$$EU(V_T^{c,\phi} - B),$$

for some (deterministic) $U$. Our arguments are applicable in this case, too:

THEOREM 2.10. *Set $U_T(\omega) = U(x - B(\omega))$. Suppose that* (NA) *and Assumption* 2.1*, as well as either Assumption* 2.3 *or Assumption* 2.5*, hold. If $B$ is bounded, then there exists $\phi^* = \phi^*(c)$ such that*

$$u(B,c) = EU(V_T^{c,\phi^*} - B) < \infty,$$

*where*

(5) $$u(B,c) := \sup_{\phi \in \Phi(U,c,B)} EU(V_T^{c,\phi} - B),$$

*and $\Phi(U,c,B)$ is the set of strategies $\phi \in \Phi$ such that the expectation $EU(V_T^{c,\phi} - B)$ is well defined.*

See the end of Section 5 for a proof.

**3. Absence of arbitrage.** Proposition 3.1 suggests that (NA) cannot be dropped in Theorem 2.7 above.

PROPOSITION 3.1. *If $U$ is strictly increasing and (NA) fails, there is no maximizer $\phi^* = \phi^*(c)$ such that $u(c) = EU(V_T^{c,\phi^*}) < +\infty$.*

PROOF. Take a strategy $\hat{\phi}$ violating (2). Then
$$EU(V_T^{c,\phi^*+\hat{\phi}}) = EU(V_T^{c,\phi^*} + V_T^{0,\hat{\phi}}) > EU(V_T^{c,\phi^*}),$$
contradicting the optimality of $\phi^*$. □

PROPOSITION 3.2. *Under (NA), the set $D_t(\omega)$ is actually a vector subspace of $\mathbb{R}^d$ almost surely.*

PROOF. This follows from Theorem 3 of [15]. It is also easy to give a direct argument. □

We will need a "quantitative" characterization of (NA), see Proposition 3.3 below. This statement is implicit in [15], but it does not follow from the arguments there. Compare also to Lemma 2.6 of [26]. Define
$$\tilde{\Xi}_t := \{\xi \in \Xi_t : |\xi(\omega)| = 1 \text{ on } \{D_{t+1} \neq \{0\}\}, \xi(\omega) \in D_{t+1}(\omega) \text{ a.s.}\}.$$

PROPOSITION 3.3. *(NA) implies the existence of $\mathcal{F}_t$-measurable random variables $\beta_t, \kappa_t > 0$ satisfying*

(6)　　$\forall p \in \tilde{\Xi}_t \quad P(\langle p, \Delta S_{t+1}\rangle < -\beta_t | \mathcal{F}_t) \geq \kappa_t \quad \text{on } \{D_{t+1} \neq \{0\}\}$

*almost surely, for all $0 \leq t \leq T-1$.*

PROOF. We may and will suppose $D_{t+1} \neq \{0\}$ a.s. Fix $t$ and a sequence $\delta_n \searrow 0$. Define
$$A_n := \left\{\omega : \operatorname*{ess\,inf}_{p \in \tilde{\Xi}_t} P(\langle p, \Delta S_{t+1}\rangle < -\delta_n | \mathcal{F}_t) = 0\right\}.$$
The essential infimum is actually attained by some $p_n^* \in \tilde{\Xi}_t$. Indeed, take $p_n^k \in \tilde{\Xi}_t$ such that
$$\lim_k \downarrow P(\langle p_n^k, \Delta S_{t+1}\rangle < -\delta_n | \mathcal{F}_t) = \operatorname*{ess\,inf}_{p \in \tilde{\Xi}_t} P(\langle p, \Delta S_{t+1}\rangle < -\delta_n | \mathcal{F}_t),$$
apply Lemma A.2 to obtain a random subsequence $\tilde{p}_n^k$ converging to some $p_n^*$. Define
$$B_k := \{\langle \tilde{p}_n^k, \Delta S_{t+1}\rangle < -\delta_n\}, \qquad B := \{\langle p_n^*, \Delta S_{t+1}\rangle < -\delta_n\},$$

and check that $B \subset \liminf_k B_k$, so $\liminf_k I_{B_k}(\omega) = I_{\liminf_k B_k}(\omega)$, and the Fatou lemma guarantees that

$$P(\langle p_n^*, \Delta S_{t+1}\rangle < -\delta_n | \mathcal{F}_t) \leq \lim_k P(\langle \tilde{p}_n^k, \Delta S_{t+1}\rangle < -\delta_n | \mathcal{F}_t),$$

so $p_n^*$ attains the essential infimum.

Clearly, $A_{n+1} \subset A_n$, set

$$A := \bigcap_{n=1}^{\infty} A_n.$$

We shall show $P(A) = 0$. If this were not the case, one would have a random subsequence $\tilde{p}_n^*$ of $p_n^*$ converging to some $\tilde{p}$. A Fatou lemma argument as above shows

$$P(\langle \tilde{p}, \Delta S_{t+1}\rangle < 0 | \mathcal{F}_t) \leq \liminf_{n \to \infty} P(\langle \tilde{p}_n^*, \Delta S_{t+1}\rangle < -\delta_n | \mathcal{F}_t) = 0$$

on $A$, so, necessarily,

$$P(\langle \tilde{p} I_A, \Delta S_{t+1}\rangle \geq 0 | \mathcal{F}_t) = 1,$$

hence, (NA) implies that

$$P(\langle \tilde{p} I_A, \Delta S_{t+1}\rangle = 0 | \mathcal{F}_t) = 1,$$

which contradicts $\tilde{p} \in D_{t+1}$, thus, indeed, $P(A) = 0$ must hold.

Define

$$\beta_t := \sum_{n=1}^{\infty} \delta_n I_{A_n^C \setminus A_{n-1}^C} \quad \text{with } A_0^C := \varnothing.$$

This is an almost everywhere positive function by $P(A) = 0$ and it is easy to see that

$$\forall p \in \tilde{\Xi}_t \quad P(\langle p, \Delta S_{t+1}\rangle < -\beta_t | \mathcal{F}_t) > 0 \quad \text{a.s.} \qquad \square$$

The condition in Proposition 3.3 is actually equivalent to (NA).

**4. Optimal strategy for the one-step case.** Let $V(x, \omega)$ be a function from $\mathbb{R} \times \Omega$ to $\mathbb{R}$ such that for almost all $\omega$, $V(\cdot, \omega)$ is a nondecreasing (finite-valued) concave function and $V(x, \cdot)$ is $\mathcal{F}$-measurable for any fixed $x$. Let $\mathcal{H} \subset \mathcal{F}$ be a $\sigma$-algebra containing $P$-zero sets. Let $Y$ be a $d$-dimensional random variable. Denote by $\Xi$ the family of $\mathcal{H}$-measurable $d$-dimensional random variables. Introduce

$$\tilde{\Xi} := \{\xi \in \Xi : |\xi(\omega)| = 1, \text{ on } \{D \neq \{0\}\}, \xi(\omega) \in D(\omega) \text{ a.s.}\},$$

here $D$ denotes the smallest affine subspace containing the support of the conditional distribution of $Y$ with respect to $\mathcal{H}$ (see the Appendix). We suppose that $D$ is actually a vector subspace a.s. and that

(7) $\quad\quad \forall p \in \tilde{\Xi} \quad\quad P(\langle p, Y \rangle < -\delta | \mathcal{H}) \geq \kappa \quad\quad \text{on } \{D \neq \{0\}\},$

with some $\mathcal{H}$-measurable random variables $\kappa, \delta > 0$.

This setting will be applied in Section 5 with the choice $\mathcal{H} = \mathcal{F}_{t-1}, D = D_t, Y = \Delta S_t$; $V(x)$ will be the maximal conditional expected utility from capital $x$ if trading begins at time $t$.

Assume that, for all $x \in \mathbb{R}$,

(8) $\quad\quad\quad\quad\quad\quad\quad\quad EV(x) > -\infty,$

and that

(9) $\quad\quad\quad\quad \underset{\xi \in \Xi}{\operatorname{ess\,sup}} E(V(x + \langle \xi, Y \rangle) | \mathcal{H}) < \infty \quad\quad \text{a.s.}$

Finally, suppose that almost surely, for all $x \in \mathbb{R}$ and $\lambda \geq 1$ both,

(10) $\quad\quad V(\lambda x) \leq \lambda V(x) + C\lambda^\gamma, \quad\quad V(\lambda x) \leq \lambda^\gamma V(x) + C\lambda^\gamma$

hold for some constants $C > 0$ and $0 < \gamma < 1$. The first inequality will be used for negative, the second for positive values of $V(x)$.

REMARK 4.1. We may interpret $V$ as a random element taking values in a Lusin space (see III. 16 of [6]): one can identify the set $\mathcal{V}$ of nondecreasing concave functions $\mathbb{R} \to \mathbb{R}$ with a Borel-subset of $\mathbb{R}^\mathbb{N}$; if $U$ is a nondecreasing concave function, then let the corresponding element of $\mathbb{R}^\mathbb{N}$ be

$$(U(q_1), U(q_2), \ldots),$$

where $(q_n)_{n \in \mathbb{N}}$ is a fixed enumeration of $\mathbb{Q}$. We leave the details to the reader, as we need this fact only once, in the proof of Proposition 4.6.

Now we briefly sketch the strategy for proving the existence of an optimal portfolio in the one-step case. After constructing regular versions of certain functions (Propositions 4.2 and 4.4), a sequence $\xi_n(x, \omega)$ is chosen along which the optimal expected utility of (9) is attained. We project strategies on $D$ (Proposition 4.6) and show [using (7) and (10)] that we may suppose $|\xi_n| \leq K$ for some $\mathcal{H}$-measurable $K$ (Lemma 4.8). Then a compactness argument provides the limit $\tilde{\xi}$ (Lemma 4.9), which turns out to be an optimal strategy.

PROPOSITION 4.2. *Let $\xi \in \Xi$ be fixed. There exists a version of*

$$x \to E(V(x + \langle \xi, Y \rangle) | \mathcal{H}),$$

*such that it is a nondecreasing upper semicontinuous concave function (perhaps taking the value $-\infty$), for almost all $\omega$.*

PROOF. We fix a version of $F(q,\omega) := E(V(q+\langle\xi,Y\rangle)|\mathcal{H})$ for $q \in \mathbb{Q}$. The following inequalities hold almost surely for any pairs $q_1 \leq q_2$ of rational numbers:

$$F(q_1) \leq F(q_2), \qquad F\left(\frac{q_1+q_2}{2}\right) \geq \frac{F(q_1)+F(q_2)}{2}.$$

Let us fix a $P$-zero set $N$ such that outside this set, the above inequalities hold. Extend $F(\cdot,\omega)$ on the real line for each $\omega \in \Omega \setminus N$ as an upper semi-continuous concave function (taking possibly the value $-\infty$). Fix $x \in \mathbb{R}$ and rationals $q_n \searrow x$. The monotone convergence theorem yields

$$F(x) = \lim_n F(q_n) = \lim_n E(V(q_n + \langle\xi,Y\rangle)|\mathcal{H}) = E(V(x+\langle\xi,Y\rangle)|\mathcal{H}),$$

showing that $F$ is, indeed, as required. □

REMARK 4.3. It is actually possible to prove the existence of a *continuous* version (taking possibly the value $-\infty$). It is also clear that the version constructed above is almost surely continuous for all $x \geq y$ if $y$ is such that

$$E(V(y+\langle\xi,Y\rangle)|\mathcal{H}) > -\infty \qquad \text{a.s.}$$

PROPOSITION 4.4. *There is a function $G : \Omega \times \mathbb{R} \to \mathbb{R}$ which is a version of*

$$\operatorname*{ess\,sup}_{\xi \in \Xi} E(V(x+\langle\xi,Y\rangle)|\mathcal{H})$$

*for each fixed $x \in \mathbb{R}$ and which is a nondecreasing finite-valued concave (a fortiori continuous) function for almost all $\omega$.*

PROOF. As in the previous proof, we construct a version $G(q,\omega)$ of the ess sup for $q \in \mathbb{Q}$ and extend it on $\mathbb{R}$ as a function which is increasing, concave and finite-valued [by (8) and (9)], hence, continuous. Fix $x \in \mathbb{R}$ and a sequence of rationals $q_n \nearrow x$. Monotone convergence shows that

$$G(x) = \lim_n \uparrow G(q_n) = \lim_n \operatorname*{ess\,sup}_{\xi \in \Xi} E(V(q_n+\langle\xi,Y\rangle)|\mathcal{H})$$

$$= \operatorname*{ess\,sup}_{\xi \in \Xi} E(V(x+\langle\xi,Y\rangle)|\mathcal{H}),$$

the proposition is proved. □

We construct a sequence of strategies converging to the optimal value for all $x \in \mathbb{R}$.

LEMMA 4.5. *There exist $\mathcal{B}(\mathbb{R}) \otimes \mathcal{H}$-measurable functions $\xi_n(x, \omega)$ and suitable versions $G_n(x)$ of*

$$E(V(x + \langle \xi_n(x), Y \rangle)|\mathcal{H}),$$

*such that outside a fixed P-zero set, we have, for all $x \in \mathbb{R}$,*

(11) $$\lim_{n \to \infty} E(V(x + \langle \xi_n(x), Y \rangle)|\mathcal{H}) = G(x),$$

*where $G(x)$ is the regular version of $\operatorname{ess\,sup}_{\xi \in \Xi} E(V(x + \langle \xi, Y \rangle)|\mathcal{H})$ figuring in Proposition 4.4. The limit is attained in a nondecreasing way.*

PROOF. It suffices to prove this for $x \in [0, 1)$; in an analogous way, we get sequences $\xi_n$ for all the intervals $[n, n+1), n \in \mathbb{Z}$ and then by "pasting together," we finally get an approximation all along the real line.

Fix a version $G(\cdot, \omega)$ of the essential supremum given by Proposition 4.4. First let us notice that, for fixed $x \in \mathbb{R}$, the family of functions

(12) $$E(V(x + \langle \xi, Y \rangle)|\mathcal{H}), \qquad \xi \in \Xi$$

is directed upwards, so there is a sequence $\eta_n(x) \in \Xi$ such that

$$\lim_{n \to \infty} \uparrow E(V(x + \langle \eta_n(x), Y \rangle)|\mathcal{H}) = \operatorname*{ess\,sup}_{\xi \in \Xi} E(V(x + \langle \xi, Y \rangle)|\mathcal{H}),$$

almost surely. Let us fix such a sequence for each dyadic rational $q \in [0, 1)$. Now set

$$\xi_0(x, \omega) := 0.$$

Suppose that $\xi_0, \ldots, \xi_{n-1}$ have been defined, as well as $\xi_n(x, \omega)$ for $0 \leq x < k/2^n$ for some $0 \leq k \leq 2^n - 1$. For $k = 0$, we set $\xi_n(x, \omega) := \kappa_n^0, x \in [0, 1/2^n)$, where $\kappa_n^0$ is chosen such that

$$E(V(\langle \kappa_n^0, Y \rangle)|\mathcal{H}) \geq E(V(\langle \xi_{n-1}(0), Y \rangle)|\mathcal{H}) \vee E(V(\langle \eta_n(0), Y \rangle)|\mathcal{H}).$$

For $k \geq 1$, we set

$$\xi_n(x, \omega) := \kappa_n^k(\omega), \qquad x \in \left[\frac{k}{2^n}, \frac{k+1}{2^n}\right),$$

where $\kappa_n^k \in \Xi$ is chosen in such a way that

(13) $$\begin{aligned} E(V(k/2^n &+ \langle \kappa_n^k, Y \rangle)|\mathcal{H}) \\ &\geq E\left(V\left(\frac{k}{2^n} + \left\langle \xi_n\left(\frac{k-1}{2^n}\right), Y \right\rangle\right)\Big|\mathcal{H}\right) \\ &\vee E\left(V\left(\frac{k}{2^n} + \left\langle \eta_n\left(\frac{k}{2^n}\right), Y \right\rangle\right)\Big|\mathcal{H}\right) \\ &\vee E\left(V\left(\frac{k}{2^n} + \left\langle \xi_{n-1}\left(\frac{k}{2^n}\right), Y \right\rangle\right)\Big|\mathcal{H}\right), \end{aligned}$$

almost everywhere. This is possible, as the family (12) is directed upwards.

Using Proposition 4.2 and Remark 4.3, take versions of the conditional expectations

$$G_n(x,\omega) := E(V(x + \langle \xi_n(x), Y \rangle)|\mathcal{H})(\omega),$$

which are nondecreasing, concave and finite-valued on intervals of the form $[k/2^n, (k+1)/2^n), 0 \leq k \leq 2^n - 1$. Proposition 4.4 and (13) show that there is a $P$-null set $N$ such that, outside this set, $G(\cdot)$ is continuous, the functions $G_n(x)$ are nondecreasing in $x$ and continuous on subintervals of the form $[k/2^n, (k+1)/2^n), 0 \leq k \leq 2^n - 1$, for $n \in \mathbb{N}$. By the definitions of $\eta_n(x)$ and $\xi_n(x)$, we see immediately that (outside another $P$-zero set $N'$), for all dyadic rationals $q$,

$$G(q) = \lim_{n \to \infty} \uparrow E(V(q + \langle \xi_n(q), Y \rangle)|\mathcal{H}) = \lim_{n \to \infty} \uparrow G_n(q).$$

Consequently, outside $N \cup N'$, the sequence $G_n(x, \omega)$ is nondecreasing in $n$, for all $x \in [0, 1)$. For any $x \in \mathbb{R}$ and dyadic rationals $q_1 < x < q_2$,

$$G_n(q_1) \leq G_n(x) \leq G_n(q_2)$$

outside $N$, so, necessarily,

$$G(q_1) \leq \liminf_n G_n(x) \leq \limsup_n G_n(x) \leq G(q_2),$$

outside $N \cup N'$. The function $G$ being continuous at $x$, we get convergence at each point $x \in [0, 1)$. □

PROPOSITION 4.6. *Let $\xi \in \Xi$. Then $\hat{\xi} \in \Xi$, where $\hat{\xi}(\omega)$ is defined as the orthogonal projection of $\xi(\omega)$ on the subspace $D(\omega)$, for all $\omega$. Furthermore,*

$$E(V(x + \langle \hat{\xi}, Y \rangle)|\mathcal{H}) = E(V(x + \langle \xi, Y \rangle)|\mathcal{H})$$

*holds almost everywhere for each $x$.*

PROOF. It is a standard exercise with the measurable selection theorem (III. 44 of [6]) to show that there exist $\mathbb{R}^d$-valued random variables $\sigma_i(\omega)$, $1 \leq i \leq d$, which almost surely span $D(\omega)$. Define the random set

$$\{(\omega, x) : x \in D(\omega), \langle \xi(\omega) - x, \sigma_i(\omega) \rangle = 0, 1 \leq i \leq d\}.$$

Almost surely it consists of one point, $\hat{\xi}(\omega)$, hence, it is the graph of a function which is measurable, again by III. 44 of [6].

We consider the random element $(V, Y) \in \mathcal{V} \times \mathbb{R}^d$ (see Remark 4.1) and denote its regular conditional probability with respect to $\mathcal{H}$ by $R(dv, dy, \omega)$ (see page 36 of [12] for an existence proof). Let its $y$-marginal be denoted by $X(dy, \omega)$. We fix any $\omega \in \Omega$ such that $R(\cdot, \cdot, \omega), X(\cdot, \omega)$ are measures. By the measure decomposition theorem of Dellacherie and Meyer ([6], III. 70–73), we have that

$$R(dv, dy, \omega) = Q(dv, y, \omega)X(dy, \omega),$$

for a suitable stochastic kernel $Q$. We have

$$E(V(x + \langle \hat{\xi}, Y \rangle)|\mathcal{H}) = \int_{\mathbb{R}^d} \int_{\mathcal{V}} v(x + \langle \hat{\xi}, y \rangle)Q(dv, y, \omega)X(dy, \omega)$$

and

$$E(V(x + \langle \xi, Y \rangle)|\mathcal{H}) = \int_{\mathbb{R}^d} \int_{\mathcal{V}} v(x + \langle \xi, y \rangle)Q(dv, y, \omega)X(dy, \omega).$$

The integrands differ only on the set $B \times \mathcal{V}$, where $B := \{y : \langle y, \xi \rangle \neq \langle y, \hat{\xi} \rangle\}$. By the definition of $D$, $X(B, \omega) = 0$, hence, the two integrals above are equal, which is just the statement of the proposition. $\square$

LEMMA 4.7.

$$\liminf_{N \to \infty} \operatorname*{ess\,inf}_{p \in \tilde{\Xi}} P(V(-N) < -1, \langle p, Y \rangle < -\delta | \mathcal{H})$$
$$\geq \operatorname*{ess\,inf}_{p \in \tilde{\Xi}} P(\langle p, Y \rangle < -\delta | \mathcal{H}).$$

PROOF. Clearly, $V(-N) \to -\infty$ almost surely as $N \to \infty$. The essential infima are attained by some $p(N)$: this can be shown just like in Proposition 3.3. So it suffices to prove

$$\liminf_{N \to \infty} P(V(-N) < -1, \langle p(N), Y \rangle < -\delta | \mathcal{H}) \geq \operatorname*{ess\,inf}_{p \in \tilde{\Xi}} P(\langle p, Y \rangle < -\delta | \mathcal{H}),$$

which follows again by taking a convergent random subsequence and the Fatou lemma. $\square$

LEMMA 4.8. *Let us fix $x_0, x_1 \in \mathbb{R}$, $x_0 < x_1$. There exists an $\mathcal{H}$-measurable random variable $K = K(x_0, x_1) > 0$ such that, for any $\xi \in \Xi$ satisfying $\xi \in D$ and $|\xi| \geq K$ a.s., we have almost surely*

$$\forall x_0 \leq x \leq x_1 \quad E(V(x + \langle \xi, Y \rangle) - V(x)|\mathcal{H}) \leq 0.$$

PROOF. Take $\xi \in \Xi, |\xi| \geq 1$ and fix a version of

$$E(V(x + \langle \xi, Y \rangle)|\mathcal{H}),$$

as given by Proposition 4.2. By (10), we have the following estimation for any $x_0 \leq x \leq x_1$:

$$V(x + \langle \xi, Y \rangle) = V^+\left(x + \left\langle \frac{\xi}{|\xi|}, Y \right\rangle |\xi| \right) - V^-\left(x + \left\langle \frac{\xi}{|\xi|}, Y \right\rangle |\xi| \right)$$

$$\leq |\xi|^\gamma V^+\left(\frac{x_1}{|\xi|} + \left\langle \frac{\xi}{|\xi|}, Y \right\rangle \right)$$

$$+ 2C|\xi|^\gamma - |\xi|^{(1+\gamma)/2} V^-\left(\frac{x_1}{|\xi|^{(1+\gamma)/2}} + \left\langle \frac{\xi}{|\xi|}, Y \right\rangle |\xi|^{(1-\gamma)/2} \right).$$

Now observe that Lemma 4.7 and (7) entail that there is an $\mathcal{H}$-measurable random variable $N_0 > 0$ such that

(14) $$\operatorname*{ess\,inf}_{p \in \tilde{\Xi}} P(V(-N_0) < -1, \langle p, Y \rangle < -\delta | \mathcal{H})(\omega) \geq \kappa/2.$$

Then

$$-E\left(V^-\left(\frac{x_1}{|\xi|^{(1+\gamma)/2}} + \left\langle \frac{\xi}{|\xi|}, Y \right\rangle |\xi|^{(1-\gamma)/2} \right) \Big| \mathcal{H}\right) \leq -E(I_B | \mathcal{H}),$$

where

$$B := \left\{ \left\langle \frac{\xi}{|\xi|}, Y \right\rangle < -\delta, V(-N_0) < -1, \frac{x_1}{|\xi|^{(1+\gamma)/2}} - |\xi|^{(1-\gamma)/2}\delta < -N_0 \right\}.$$

Putting together our estimations so far,

(15) $$E(V(x + \langle \xi, Y \rangle) | \mathcal{H}) \leq |\xi|^\gamma E\left(V^+\left(x_1 + \left\langle \frac{\xi}{|\xi|}, Y \right\rangle \right) \Big| \mathcal{H}\right)$$

(16) $$+ 2C|\xi|^\gamma - |\xi|^{(1+\gamma)/2}\kappa/2,$$

as soon as

(17) $$\frac{x_1}{|\xi|^{(1+\gamma)/2}} - |\xi|^{(1-\gamma)/2}\delta < -N_0.$$

Let $K_0(\omega) > 0$ be an $\mathcal{H}$-measurable random variable such that (17) is true for $|\xi| \geq K_0$. We shall show in a minute that the first term on the right-hand side of (15) is smaller than $L|\xi|^\gamma$ for some $\mathcal{H}$-measurable $L(\omega) > 0$, so this right-hand side is smaller than $E(V(x_0)|\mathcal{H})$, provided that $|\xi| \geq K_1$, for some $\mathcal{H}$-measurable $K_1(\omega) > 0$. Finally, we may conclude that if $K := K_0 \vee K_1 \vee 1$, then the statement of Lemma 4.8 holds.

It remains to estimate the first term of the right-hand side of (15). Introduce the following vectors for $i \in W := \{-1, +1\}^d$:

$$\theta_i^j := i(j), \qquad j = 1, \ldots, d.$$

It is easy to see that

$$V^+(x_1 + \langle p, Y \rangle) \leq \max_{i \in W} V^+(x_1 + \langle \theta_i, Y \rangle),$$

for any $p \in \mathbb{R}^d, |p| \leq 1$. Hence, the term in consideration is smaller than

$$(18) \qquad |\xi|^\gamma \sum_{i \in W} E(V^+(x_1 + \langle \theta_i, Y \rangle)|\mathcal{H}) =: |\xi|^\gamma L(\omega),$$

and $L$ is finite by (9). □

LEMMA 4.9.  *There exists a $\mathcal{B}(\mathbb{R}) \otimes \mathcal{H}$-measurable function $\tilde{\xi}(x, \omega)$ such that, for all $x \in \mathbb{R}$,*

$$E(V(x + \langle \tilde{\xi}(x), Y \rangle)|\mathcal{H}) = \operatorname*{ess\,sup}_{\xi \in \Xi} E(V(x + \langle \xi, Y \rangle)|\mathcal{H}).$$

PROOF.  It suffices to prove this, for example, $x \in [0, 1)$, then one can "paste together" the optimal strategy for $x \in \mathbb{R}$. We take an approximating sequence $\xi_n$ as provided by Lemma 4.5, then change to the projections $\hat{\xi}_n$ figuring in Proposition 4.6. Using Proposition 4.2 and the structure, of the approximating sequence, one can see that there are suitable versions of

$$E(V(x + \langle \hat{\xi}_n(x), Y \rangle)|\mathcal{H}),$$

such that almost surely,

$$\forall x \in [0, 1) \qquad E(V(x + \langle \hat{\xi}_n(x), Y \rangle)|\mathcal{H}) \to G(x), \qquad n \to \infty.$$

Then take $x_0 := 0, x_1 := 1$ and truncate $\hat{\xi}_n$: the strategies

$$\eta_n := \hat{\xi}_n I_{\{|\hat{\xi}_n| \leq K(x_0, x_1)\}}$$

do at least as well as the original sequence since by Lemma 4.8 (and using suitable versions of the conditional expectations) almost surely,

$$\forall x \in [0, 1) \qquad E(V(x + \langle \hat{\xi}_n(x), Y \rangle)|\mathcal{H}) \leq E(V(x + \langle \eta_n(x), Y \rangle)|\mathcal{H}).$$

Again, for suitable versions of the conditional expectations, we almost surely have

$$\forall x \in [0, 1) \qquad E(V(x + \langle \eta_n(x), Y \rangle)|\mathcal{H}) \to G(x), \qquad n \to \infty.$$

Now use Lemma A.2 to find a random subsequence $\tilde{\eta}_k = \eta_{n_k}$ of $\eta_n$ converging to some $\tilde{\xi}$. Apply the Fatou lemma (we shall justify its use in a while):

$$E(V(x + \langle \tilde{\xi}(x), Y \rangle)|\mathcal{H}) \geq \limsup_k E(V(x + \langle \tilde{\eta}_k(x), Y \rangle)|\mathcal{H}) \qquad \text{a.s.,}$$

for each fixed $x$. By construction, almost surely

$$\forall x \qquad E(V(x + \langle \tilde{\eta}_k(x), Y \rangle)|\mathcal{H}) \geq E(V(x + \langle \xi_{n_k}(x), Y \rangle)|\mathcal{H}),$$

so the definition of essential supremum and the construction imply that, for each fixed $x \in \mathbb{R}$,
$$E(V(x + \langle \tilde{\xi}(x), Y \rangle)|\mathcal{H}) = G(x),$$
that is, $G(x)$ is a version of the conditional expectation.

It remains to check that one is allowed to use the Fatou lemma. This can be shown as in Lemma 4.8: take the $\theta_i, i \in W$ defined there and estimate as follows:
$$V^+(x + \langle \tilde{\eta}_n, Y \rangle) \leq \max_{i \in W} V^+(x + K\langle \theta_i, Y \rangle) \leq \sum_{i \in W} V^+(x + K\langle \theta_i, Y \rangle),$$
for each $n$, hence,
$$E\left(\max_n V^+(x + \langle \tilde{\eta}_n, Y \rangle)|\mathcal{H}\right) \leq \sum_{i \in W} E(V^+(x + K\langle \theta_i, Y \rangle)|\mathcal{H}) < \infty,$$
due to (9). □

PROPOSITION 4.10. *The $\tilde{\xi}$ constructed in the proof of Lemma 4.9 satisfies*
$$G(H) = E(V(H + \langle \tilde{\xi}(H), Y \rangle)|\mathcal{H}) = \operatorname*{ess\,sup}_{\xi \in \Xi} E(V(H + \langle \xi, Y \rangle)|\mathcal{H}) \quad a.s.,$$
*for any $\mathcal{H}$-measurable $\mathbb{R}$-valued random variable $H$; here $G$ is the function constructed in Proposition 4.4.*

PROOF. One may suppose, for example, $H \in [0, 1)$. Fix $n \in \mathbb{N}$. Clearly,
$$(19) \qquad G_n(H) = E(V(H + \langle \xi_n(H), Y \rangle)|\mathcal{H}) \quad \text{a.s.,}$$
for step functions $H$, see the proof of Lemma 4.5. For general $H$, take step function approximations $H_l \searrow H, l \to \infty$ such that $H_l \in [k/2^n, (k+1)/2^n)$ on the set $\{\omega : H(\omega) \in [k/2^n, (k+1)/2^n)\}$, for all $0 \leq k \leq 2^n - 1$. The strategies $\xi_n$ are piecewise constant, $G_n(H_l) \to G_n(H), l \to \infty$ by piecewise continuity, so monotone convergence implies (19) for general $H$. Now the proof of Lemma 4.9 shows
$$G_{n_k}(H) \leq E(V(H + \langle \eta_{n_k}(H), Y \rangle)|\mathcal{H})$$
almost surely, for each $k \in \mathbb{N}$. Letting $k \to \infty$, we get, by Lemma 4.5,
$$G(H) \leq E(V(H + \langle \tilde{\xi}(H), Y \rangle)|\mathcal{H}) \quad \text{a.s.}$$

The left-hand side of the second equality in the statement of Proposition 4.10 is clearly not greater than the right-hand side, so we only need to show that for fixed $\xi \in \Xi$,
$$(20) \qquad G(H, \omega) \geq E(V(H + \langle \xi, Y \rangle)|\mathcal{H}) \quad \text{a.s.}$$
For step functions $H$, (20) is clearly true. Taking step functions $H_n \searrow H$, the left-hand side converges by path regularity of $G$, the right-hand side by monotone convergence. □

**5. Dynamic programming.** First we need an easy fact about $U$.

PROPOSITION 5.1. *Let $U$ satisfy Assumption* 2.3. *Then there is a constant $C > 0$ such that $U$ satisfies for all $x$ and all $\lambda \geq 1$ both of the following inequalities*:

$$U(\lambda x) \leq \lambda U(x) + C\lambda^\gamma, \tag{21}$$

$$U(\lambda x) \leq \lambda^\gamma U(x) + C\lambda^\gamma. \tag{22}$$

PROOF. Obviously (21) holds true for $x \geq \tilde{x}$, since $\lambda \geq \lambda^\gamma$ and $0 < \gamma < 1$. As $U$ is nondecreasing, (2.3) implies that

$$U(\lambda x) \leq U(\lambda \tilde{x}) \leq \lambda^\gamma U(\tilde{x}), \tag{23}$$

for $0 \leq x \leq \tilde{x}$, so we may set $C := U(\tilde{x})$. For $x < 0$, we have the following estimation by concavity:

$$U(\lambda x) \leq U(x) + U'(x)(\lambda - 1)x \leq U(x) + (\lambda - 1)(U(x) - U(0)) = \lambda U(x).$$

Now (22) is clear for $x > 0$ by Assumption 2.3 and (23). Finally, (22) for $x < 0$ follows from (21), since in this case $U(x) \leq 0$ and $\lambda^\gamma \leq \lambda$. □

We would like to perform a dynamic programming argument in a non-Markovian context just like Evstigneev [7]; establishing that some crucial properties of $U$ are preserved by $U_t$. In particular, the "asymptotic elasticity" conditions (21) and (22). In continuous-time models such preservation properties are studied in Lemma 3.12 of [20] and in section "Dynamic version of the utility maximisation problem" of [27].

PROPOSITION 5.2. *The functions $U_t, 0 \leq t \leq T$ have versions, which are almost surely nondecreasing, concave, and satisfy* (21) *and* (22), *as well as*

$$EU_t(x) > -\infty, \qquad x \in \mathbb{R}, 0 \leq t \leq T, \tag{24}$$

$$U_t(x) < \infty, \qquad x \in \mathbb{R}, 0 \leq t \leq T. \tag{25}$$

*There exist $\mathcal{B}(\mathbb{R}) \otimes \mathcal{F}_t$-measurable functions $\tilde{\xi}_{t+1}, 1 \leq i \leq T$ such that*

$$\forall x \in \mathbb{R} \qquad U_t(x,\omega) = E(U_{t+1}(x + \langle \tilde{\xi}_{t+1}(x), \Delta S_{t+1}\rangle)|\mathcal{F}_t). \tag{26}$$

PROOF. Going backwards from $T$ to 0, apply Lemma 4.9 with the choice $V := U_t, \mathcal{H} = \mathcal{F}_{t-1}, D := D_t, Y := \Delta S_t$. We need to verify the conditions of Section 4: $D$ is a random subspace by Propositions 3.2 and A.1; (7) follows from (6); (8) and (9) will come from (24) and (25); (10) is a consequence of (21) and (22); (25) follows from Assumption 2.1. We will check (21), (22) and (24) in a little while. Denote the resulting $\tilde{\xi}$ of Lemma 4.9 by $\tilde{\xi}_t$, $1 \leq t \leq T$, it satisfies (26).

Good versions exist by Proposition 4.4. For $t = T$, (22) holds because of Proposition 5.1 and for $t < T$, by

$$U_{t-1}(\lambda x) = E(U_t(\lambda x + \langle \tilde{\xi}_t(\lambda x), \Delta S_t \rangle) | \mathcal{F}_{t-1})$$
$$\leq \lambda^\gamma (E(U_t(x + \langle \tilde{\xi}_t(\lambda x)/\lambda, \Delta S_t \rangle) | \mathcal{F}_{t-1}) + C) \leq \lambda^\gamma (U_{t-1}(x) + C).$$

We get (21) in the same way. It remains to establish (24): the statement is true, since

(27) $$U_t(x) \geq E(U_{t+1}(x) | \mathcal{F}_t) \geq \cdots \geq U_T(x),$$

and the latter is deterministic. □

Now set $\phi_1^* := \tilde{\xi}_1(c)$ and define inductively

$$\phi_{t+1}^* := \tilde{\xi}_{t+1}\left(c + \sum_{j=1}^t \langle \phi_j^*, \Delta S_j \rangle\right), \qquad 1 \leq t \leq T-1.$$

Joint measurability of $\tilde{\xi}_t$ assures that $\phi^*$ is a predictable process with respect to the given filtration.

PROPOSITION 5.3. *For any strategy* $\phi \in \Phi(U, x)$,

(28) $$E(U(V_T^{c,\phi}) | \mathcal{F}_0) \leq E(U(V_T^{c,\phi^*}) | \mathcal{F}_0) = U_0(c).$$

PROOF. Remembering $U_T = U$ and using Proposition 4.10, we may rewrite the right-hand side of (28) as follows:

$$E(U_T(V_T^{c,\phi^*}) | \mathcal{F}_0) = E(E(U_T(V_{T-1}^{c,\phi^*} + \langle \phi_T^*, \Delta S_T \rangle) | \mathcal{F}_{T-1}) | \mathcal{F}_0)$$
$$= E(U_{T-1}(V_{T-1}^{c,\phi^*}) | \mathcal{F}_0).$$

Continuing the procedure, we finally arrive at

(29) $$E(U(V_T^{c,\phi^*}) | \mathcal{F}_0) = E(U_1(V_1^{c,\phi^*}) | \mathcal{F}_0)$$
$$= E(U_1(c + \langle \phi_1^*, \Delta S_1 \rangle) | \mathcal{F}_0) = U_0(c).$$

We remark that all conditional expectations below exist by the definition of $\Phi(U, x)$. By the definition of $U_{T-1}$,

$$E(U_T(V_T^{c,\phi}) | \mathcal{F}_{T-1}) = E(U_T(V_{T-1}^{c,\phi} + \langle \phi_T, \Delta S_T \rangle) | \mathcal{F}_{T-1}) \leq U_{T-1}(V_{T-1}^{c,\phi}).$$

Iterate the same argument and obtain

(30) $$E(U(V_T^{c,\phi}) | \mathcal{F}_0) \leq U_0(c).$$

Putting (29) and (30) together, one gets exactly (28). □

PROOF OF THEOREM 2.7 UNDER ASSUMPTION 2.3. Proposition 5.3 shows that $u(c) = EU_0(c)$ and $\phi^*(c)$ is an optimal strategy. □

PROOF OF THEOREM 2.7 UNDER ASSUMPTION 2.5. Define

$$\tilde{U}(x) := U(x + \tilde{x}) - U(\tilde{x}).$$

Assumption 2.1 holds for this new function and an optimal strategy for $\tilde{U}$ furnishes one for $U$. $\tilde{U}(0) = 0$ and if $x \leq 0$, we have, by Assumption 2.5 for $\lambda \geq 1$,

$$\tilde{U}(\lambda x) \leq \lambda^{1+\alpha} U\left(x + \frac{\tilde{x}}{\lambda}\right) - U(\tilde{x})$$

$$\leq \lambda^{1+\alpha} \tilde{U}\left(x - \left(1 - \frac{1}{\lambda}\right)\tilde{x}\right).$$

Let us introduce the notation

$$\chi := -\tilde{x} \geq -\left(1 - \frac{1}{\lambda}\right)\tilde{x},$$

so we have

$$\tilde{U}(\lambda x) \leq \lambda^{1+\alpha} \tilde{U}(x + \chi),$$

for $x \leq 0$. An argument similar to that of Proposition 5.1 above and Lemma 6.3 of [20] shows that, for all $x \in \mathbb{R}$,

$$\tilde{U}(\lambda x) \leq \lambda^{1+\alpha} \tilde{U}(x + \chi),$$
$$\tilde{U}(\lambda x) \leq \lambda \tilde{U}(x + \chi).$$

Replacing Proposition 5.1 by the two inequalities above, only minor modifications are needed in the estimations of Lemma 4.8 and Proposition 5.2, otherwise the proof of Theorem 2.7 goes through for $\tilde{U}$. □

PROOF OF THEOREM 2.8. If $U$ is strictly concave, the functions $U_t$ are easily seen to be strictly concave too. Now unicity has to be proved by forward induction. Let us suppose that $\phi_i \in \Xi_0, \phi_i \in D_1$ a.s., $i = 1, 2$, such that

$$E(U_1(c + \langle \phi_i, \Delta S_1 \rangle)) = u(c)$$

for $i = 1, 2$ with a given fixed $c$. Then for the strategy $\phi_3 := (\phi_1 + \phi_2)/2$, we obtain

$$EU_1(c + \langle \phi_3, \Delta S_1 \rangle) \geq \frac{EU_1(c + \langle \phi_1, \Delta S_1 \rangle) + EU_1(c + \langle \phi_2, \Delta S_1 \rangle)}{2} \geq u(c),$$

by concavity. By the definition of $u(c)$, only equality is possible, hence, strict concavity implies that, necessarily,

$$\langle \phi_3, \Delta S_1 \rangle = \langle \phi_2, \Delta S_1 \rangle = \langle \phi_1, \Delta S_1 \rangle \quad \text{a.s.}$$

But this implies that actually $\phi_1 - \phi_2 \in D_1^\perp$, which is only possible if

$$\phi_1 - \phi_2 = 0 \quad \text{a.s.,}$$

as required. The induction step is identical, one has to apply the induction hypothesis and consider

$$EU_t(V_{t-1}^{c,\phi^*} + \langle \phi_i, \Delta S_t \rangle)$$

for $\phi_i \in \Xi_{t-1}, \phi_i \in D_t$ a.s., $i = 1, 2$. The above argument shows $\phi_1 = \phi_2$ almost surely. $\square$

PROOF OF THEOREM 2.10. First suppose Assumption 2.3. Let us define

$$U_T(x, \omega) := U(x - B(\omega)),$$

and $U_t, 0 \leq t \leq T-1$ in a respective manner. Boundedness of $B$ and Assumption 2.1 imply (1). Furthermore, observe that if $\ell \in \mathbb{R}$ is such that $\ell \geq |B|$ almost surely, then by Proposition 5.1,

$$U_T(\lambda x) = U(\lambda x - B) \leq \lambda U(x - B/\lambda) + C\lambda^\gamma$$
$$\leq \lambda U(x - B + \ell) + C\lambda^\gamma = \lambda U_T(x + \ell) + C\lambda^\gamma$$

and

$$U_T(\lambda x) \leq \lambda^\gamma U_T(x + \ell) + C\lambda^\gamma,$$

for all $\lambda \geq 1$. Hence, it is easy to see that, apart from trivial changes in the estimations of Lemma 4.8 and Proposition 5.2, the arguments of Sections 4 and 5 go through. Under Assumption 2.5, we use

$$\tilde{U}(x) := U_T(x + \tilde{x} - \ell - B) - U_T(\tilde{x} - \ell - B),$$

and the proof of Theorem 2.7 under this assumption (see above) gives

$$\tilde{U}(\lambda x) \leq \lambda^{1+\alpha} \tilde{U}(x - \tilde{x} + 2\ell),$$
$$\tilde{U}(\lambda x) \leq \lambda \tilde{U}(x - \tilde{x} + 2\ell),$$

the same arguments work almost without modifications. $\square$

**6. Utility based pricing.** We are looking for equivalent martingale measures, that is, $Q \sim P$ such that $S$ is a $Q$-martingale, using utility considerations. This approach has already been pursued by, for example, Davis [5], Frittelli [8] and Hugonnier, Kramkov and Schachermayer [14], and goes back to principles of economic theory, see [13].

Recently Schäl [[29]–[31]] have investigated this method in the discrete-time context. A natural candidate for an equivalent martingale measure is

$$\text{(31)} \qquad \frac{dQ}{dP} := \frac{U'(V_T^{c,\phi^*})}{EU'(V_T^{c,\phi^*})},$$

where $\phi^*$ is the optimal strategy for initial capital $c$.

REMARK 6.1. Unfortunately, in several cases, formula (31) fails to provide a martingale measure, even for a "nice" $U$. Let us fix, for example, $U(x) := 1 - e^{-x}, T := 1, c := 0$ and let $\Delta S_1$ be symmetric and not integrable. Then

$$Ee^{-\phi \Delta S_1} = \infty, \qquad \phi \in \mathbb{R} \setminus \{0\},$$

the optimal strategy is given by $\phi_1^* := 0$, the corresponding $Q$ is $P$ itself; but $P$ is certainly not a martingale measure!

THEOREM 6.2. *Let us suppose that $U$ satisfies either Assumption 2.3 or Assumption 2.5 and that it is strictly increasing and continuously differentiable. Furthermore, assume that $S$ satisfies (2) and is bounded. If Assumption 2.1 is met, then (31) defines an equivalent martingale measure.*

Before proving this theorem, we need several auxiliary assertions. First notice that, by Theorem 2.7, there is an optimal $\phi^*$.

PROPOSITION 6.3. *Let $f : \mathbb{R}^d \to \mathbb{R}$ a concave function. Then*

$$|f(s) - f(s')| \leq |s - s'| \left\{ \left| f(s) - f\left(s + \frac{s' - s}{|s - s'|}\right) \right| \vee \left| f(s') - f\left(s' + \frac{s - s'}{|s - s'|}\right) \right| \right\},$$

*for $|s - s'| \leq 1$.*

PROOF. We see that, for any $\lambda \in [0, 1]$,

$$f(s) - f(s') = f(s) - f\left(\lambda s + \frac{1 - \lambda}{1 - \lambda}(s' - \lambda s)\right),$$

so concavity implies that the left-hand side is not greater than

$$f(s) - \lambda f(s) - (1 - \lambda) f\left(\frac{s' - \lambda s}{1 - \lambda}\right) = (1 - \lambda)\left(f(s) - f\left(\frac{s' - \lambda s}{1 - \lambda}\right)\right).$$

Let $|s - s'| \leq 1$ and let us set $\lambda := 1 - |s - s'|$. We obtain

$$f(s) - f(s') \leq |s - s'|\left(f(s) - f\left(s + \frac{s' - s}{|s - s'|}\right)\right).$$

The statement of Proposition 6.3 now follows by interchanging the roles of $s$ and $s'$. □

From now on, we use notation and assumptions of Section 4, as well as hypotheses of Theorem 6.2. In addition, suppose that for almost all $\omega$, $V(\cdot, \omega)$ is continuously differentiable, $Y$ is bounded and

(32) $$-\infty < EV(x + \langle y, Y \rangle) < \infty,$$

for all $x \in \mathbb{R}, y \in \mathbb{R}^d$.

PROPOSITION 6.4. *The function*

$$(x, y) \to E(V(x + \langle y, Y \rangle)|\mathcal{H})$$

*has a version $H(x, y, \omega)$ which is continuously differentiable in $(x, y) \in \mathbb{R}^{d+1}$,*

$$\partial_i H(x, y, \omega) = E(V'(x + \langle y, Y \rangle)Y^i|\mathcal{H}), \qquad 1 \leq i \leq d,$$

*where $\partial_i$ is the derivative with respect to $y^i$,*

(33) $$\partial_x H(x, y, \omega) = E(V'(x + \langle y, Y \rangle)|\mathcal{H}).$$

*Furthermore, for any $\xi \in \Xi$ and any $\mathcal{H}$-measurable $\mathbb{R}$-valued r.v. $X$,*

$$H(X, \xi, \omega) = E(V(X + \langle \xi, Y \rangle)|\mathcal{H}),$$
$$\partial_i H(X, \xi, \omega) = E(V'(X + \langle \xi, Y \rangle)Y^i|\mathcal{H}), \qquad 1 \leq i \leq d.$$

*So one has also*

$$E(V'(X + \langle \tilde{\xi}(X), Y \rangle)Y^i|\mathcal{H}) = 0, \qquad 1 \leq i \leq d.$$

PROOF. We confine ourselves to the case $d := 1$ and $x, y \in [0, 1]$. First apply Proposition 6.3 with the choice $s' := x + yY, s := x + (y + h)Y$ for $h \in \mathbb{R}$ such that $|hY| \leq 1$:

$$\frac{|V(x + (y + h)Y) - V(x + yY)|}{|h|}$$

(34) $$\leq |Y|\{|V(x + (y + h)Y)| + |V(x + yY)|\}$$
$$+ |Y|\{|V(x + yY + 1)| + |V(x + (y + h)Y + 1)|\}.$$

Condition (32) implies that the right-hand side of the above inequality is in $L^1$. Hence,

$$V'(x + yY)Y = \lim_{h \to 0} \frac{V(x + (y + h)Y) - V(x + hY)}{h}$$

is also in $L^1$, even $\sup_{x,y \in [0,1]} |V'(x+yY)Y|$ is in $L^1$. A similar argument works for (33).

Now apply Lemma A.3 and obtain a continuously differentiable version of $H$.

Notice that the second assertion is clear for step functions $X, \xi$ and we may also suppose $X, \xi \geq 0$. Taking arbitrary $\xi \geq 0$ and increasing step-function approximations $\xi_n \nearrow \xi$, we get

$$\partial_1 H(X, \xi_n) \to \partial_1 H(X, \xi)$$

by continuity and

$$\begin{aligned} E(V'(X+\xi_n Y)Y|\mathcal{H}) &= E(V'(X+\xi_n Y^+)Y^+|\mathcal{H}) - E(V'(X-\xi_n Y^-)Y^-|\mathcal{H}) \\ &\to E(V'(X+\xi Y)Y|\mathcal{H}) \end{aligned}$$

by monotone convergence. The above reasoning for $X_n \nearrow X$ and $\xi$ fixed completes the argument. The analogous statement for $H$ follows in a similar way.

In the present case, Lemma 4.9 gives a $P$-zero set $N$ such that

$$(35) \qquad \forall x \quad H(x, \tilde{\xi}(x), \omega) = G(x, \omega), \qquad \omega \in \Omega \setminus N,$$

where $G$ is as constructed in Proposition 4.4.

From the definition of $G$ and by continuity of $G, H$,

$$(36) \qquad \forall x, y \quad H(x, y, \omega) \leq G(x, \omega), \qquad \omega \in \Omega \setminus N',$$

outside another zero-set $N'$. Now the last assertion follows by the optimality of $\tilde{\xi}$. □

We improve on Proposition 4.4 next.

PROPOSITION 6.5. *The function*

$$x \to \operatorname*{ess\,sup}_{\xi \in \Xi} E(V(x + \langle \xi, Y \rangle)|\mathcal{H})$$

*has a version $G(x, \omega)$ which is almost surely continuously differentiable and, for any $\mathcal{H}$-measurable $\mathbb{R}$-valued random variable $X$,*

$$G(X, \omega) = \operatorname*{ess\,sup}_{\xi \in \Xi} E(V(X + \langle \xi, Y \rangle)|\mathcal{H}).$$

*Also,*

$$(37) \qquad G'(X, \omega) = E(V'(X + \langle \tilde{\xi}(X), Y \rangle)|\mathcal{H}).$$

PROOF. Take $G$ as given by Proposition 4.4. From Proposition 6.4,

$$\partial_s H(s,y,\omega) = E(V'(s + \langle y, Y \rangle)|\mathcal{H}).$$

Now we borrow a trick from Theorem 4.13 of [30]. Outside a null set, for all $x$ and for any $h \in \mathbb{N}$,

$$G(x \pm 1/h) - G(x) \geq H(x \pm 1/h, \tilde{\xi}(x)) - H(x, \tilde{\xi}(x)),$$

see the end of the proof of Proposition 6.4. Letting $h \to \infty$, we find that

$$\partial_x H(x, \tilde{\xi}(x)) \geq G'(x-) \geq G'(x+) \geq \partial_x H(x, \tilde{\xi}(x)),$$

by smoothness of $H$, so $G$ is indeed smooth almost everywhere, a similiar argument assures that one can plug $X$ into $G'$. □

PROPOSITION 6.6. *The functions $U_t, 0 \leq t \leq T$, have continuously differentiable versions which also satisfy (32). Furthermore, we have for $1 \leq i \leq d$ and for $1 \leq t \leq T$,*

$$E(U_t'(X + \langle \tilde{\xi}_t(X), \Delta S_t \rangle)\Delta S_t^i | \mathcal{F}_{t-1}) = 0,$$

*for any $\mathcal{F}_{t-1}$-measurable $X$.*

PROOF. For $t = T$, the first two assertions are clear as $S$ is bounded. The rest follows by Propositions 6.4 and 6.5 and backward induction; (32) holds true because of (8) and Assumption 2.1. □

PROOF OF THEOREM 6.2. We need to check that, for all $0 \leq t \leq T - 1$,

$$E(U'(V_T^{c,\phi^*})|\mathcal{F}_t) = U_t'(V_t^{c,\phi^*}).$$

Indeed, this follows by backward induction and (37). We also get, by an estimation like (34), that

$$U_0'(c) = E(U'(V_T^{c,\phi^*})|\mathcal{F}_0)$$

is in $L^1$, thus, Proposition 6.6 implies that the $Q$ defined by (31) is an absolutely continuous martingale measure. As $U$ is strictly increasing and concave, $U'$ never vanishes, so $Q$ is equivalent to $P$. □

**7. Ramifications.** We would like to check that Theorems 2.7 and 6.2 hold in concrete, nontrivial classes of models. Let $\mathcal{M}$ denote the set of random variables with finite moments of all orders.

PROPOSITION 7.1. *Let us suppose Assumption 2.3 (or Assumption 2.5),*

(38) $$M|x|^{-l} \leq U'(x) \leq K(|x|^k + 1)$$

*for some* $k, l, M, K \geq 0$; *and that* $U$ *is continuously differentiable and strictly increasing. Furthermore, assume that for all* $0 \leq t \leq T$, *we have* $|\Delta S_t| \in \mathcal{M}$ *and* (NA) *holds such that* $\beta_t, \kappa_t$ *of Proposition* 3.3 *satisfy* $1/\beta_t, 1/\kappa_t \in \mathcal{M}$ *for* $0 \leq t \leq T-1$ *(this applies, in particular, when* $\kappa_t = \kappa, \beta_t = \beta$ *deterministic constants, e.g., when* $S$ *has independent increments).*

*Then Assumption* 2.1 *holds; for every initial endowment* c, *there exists a strategy* $\phi^*(c)$ *such that*

$$u(c) = EU(V_T^{c,\phi^*}) < \infty,$$

*and* (31) *defines an equivalent martingale measure.*

PROOF. We shall check that the arguments in the proofs of Theorems 2.7 and 6.2 work. The main point consists in establishing a more "quantitative" version of Lemma 4.8. Suppose Assumption 2.3, the case of Assumption 2.5 is analogous.

We shall show by backward induction that, for all $0 \leq t \leq T$ and $x \in \mathbb{R}$,

(39) $$U_t(x) \leq E(U(x + (1+|x|^{\zeta_t})\rho_t)|\mathcal{F}_t),$$

(40) $$U_t(x) < \infty \quad \text{a.s.}$$

(41) $$\exists \tilde{\xi}_{t+1}(x) \in \Xi_t \quad U_t(x) = E(U_{t+1}(x + \langle \tilde{\xi}_{t+1}(x), \Delta S_{t+1}\rangle)|\mathcal{F}_t),$$

(42) $$|\tilde{\xi}_{t+1}(x)| \leq (1+|x|^{\alpha_t})\psi_t,$$

(43) $$U_t \text{ is continuously differentiable},$$

(44) $$U'_t(x) = E(U'_{t+1}(x + \langle \tilde{\xi}_{t+1}(x), \Delta S_{t+1}\rangle)|\mathcal{F}_t),$$

where $\zeta_t, \alpha_t > 0$ constants, $\rho_t, \psi_t \in \mathcal{M}$.

Suppose that the above statements are true for $t+1$ and proceed by the induction step (the case $t = T$ is trivial). Estimations of Lemma 4.8 [(14), (15), (17) and (18), in particular] show that, for any $\phi \in \Xi_t, \phi \in D_{t+1}, |\phi| \geq 1$,

$$E(U_{t+1}(x + \langle \phi, \Delta S_{t+1}\rangle)|\mathcal{F}_t) \leq |\phi|^\gamma L(x) + 2C|\phi|^\gamma - |\phi|^{(1+\gamma)/2}\kappa_t/2,$$

whenever

(45) $$\operatorname*{ess\,inf}_{q \in \tilde{\Xi}_t} P(U_{t+1}(|x| - |\phi|^{(1-\gamma)/2}\beta_t) < -1, \langle q, \Delta S_{t+1}\rangle < -\beta_t|\mathcal{F}_t) \geq \kappa_t/2,$$

and $L(x)$ is defined as

$$L(x) := \sum_{i \in W} E(U_{t+1}^+(x + \langle \theta_i, \Delta S_{t+1}\rangle)|\mathcal{F}_t).$$

Estimating $L(x)$ from above by (39) as

$$L(x) \leq 2^d E(U^+(|x| + [1 + (|x| + \sqrt{d}|\Delta S_{t+1}|)^{\zeta_{t+1}}]\rho_{t+1})|\mathcal{F}_t)$$
$$\leq 2^d U'(0) E(|x| + [1 + (|x| + \sqrt{d}|\Delta S_{t+1}|)^{\zeta_{t+1}}]\rho_{t+1}|\mathcal{F}_t),$$

we conclude that $L(x) \leq (1+|x|^w)J$, where $w > 0, J \in \mathcal{M}$. Turning our attention to (45), let us use concavity, the induction hypotheses, (39), (42) and (44), while supposing $|\phi|^{(1-\gamma)/2}\beta_t \geq |x|$:

$$U_{t+1}(|x| - |\phi|^{(1-\gamma)/2}\beta_t)$$
$$\leq U_{t+1}(0) - U'_{t+1}(0)(|\phi|^{(1-\gamma)/2}\beta_t - |x|)$$
$$\leq E(U(\rho_{t+1})|\mathcal{F}_{t+1}) - E(U'(\Psi)|\mathcal{F}_{t+1})(|\phi|^{(1-\gamma)/2}\beta_t - |x|),$$

where $\Psi \in \mathcal{M}$. Define $N_1, N_2 \in \mathcal{M}$ by

$$N_1 := \frac{4E(\Psi|\mathcal{F}_t)}{\kappa_t}, \qquad N_2 := \frac{4E(U(\rho_{t+1})|\mathcal{F}_t)}{\kappa_t}.$$

The Markov inequality then assures

$$P(\Psi \leq N_1, E(U(\rho_{t+1})|\mathcal{F}_{t+1}) \leq N_2 | \mathcal{F}_t) \geq 1 - \kappa_t/2.$$

Hence, by (6) and the above considerations,

$$\operatorname*{ess\,inf}_{q \in \tilde{\Xi}_t} P(U_{t+1}(|x| - |\phi|^{(1-\gamma)/2}\beta_t) < -1, \langle q, \Delta S_t \rangle < -\beta_t | \mathcal{F}_t)$$

$$\geq \operatorname*{ess\,inf}_{q \in \tilde{\Xi}_t} P\Big(U'(N_1)(\beta_t |\phi|^{(1-\gamma)/2} - |x|) > 1 + N_2,$$

$$E(U(\rho_{t+1})|\mathcal{F}_{t+1}) \leq N_2, \Psi \leq N_1, \langle q, \Delta S_t \rangle < -\beta_t | \mathcal{F}_t\Big) \geq \kappa_t/2,$$

provided that

(46) $$|\phi| \geq \left[\left(\frac{1+N_2}{U'(N_1)} + |x|\right)/\beta_t\right]^{2/(1-\gamma)}.$$

Note that, due to (38), condition (46) is met as soon as $|\phi| \geq (1+|x|^s)\Theta$ for some $s > 0, \Theta \in \mathcal{M}$.

Choose $\phi$ so large as to have

$$|\phi|^\gamma L + 2C|\phi|^\gamma - |\phi|^{(1+\gamma)/2}\kappa_t/2$$
$$< -\frac{K}{k+1}E([|x| + (1+|x|^{\zeta_{t+1}})\rho_{t+1}]^{k+1}|\mathcal{F}_t) \leq E(U_{t+1}(x)|\mathcal{F}_t),$$

where we used (39) and (38) in the second inequality. Again, by the bound on $L(x)$, this gives a condition on $|\phi|$ which is polynomial in $x$ and involves terms in $\mathcal{M}$, that is, the essential supremum $U_t$ is attained by portfolios satisfying (42) for appropriate $\alpha_t, \psi_t$ given by the tedious estimations above. Now it follows that if $\phi$ satisfies the bound (42) then for suitable $\zeta_t > 0, \rho_t \in \mathcal{M}$,

$$E(U_{t+1}(x + \langle \phi, \Delta S_t \rangle)|\mathcal{F}_t)$$
$$\leq E(U(x + (1 + [|x| + |\phi||\Delta S_t|]^{\zeta_{t+1}})\rho_{t+1})|\mathcal{F}_t)$$
$$\leq E(U(x + (1+|x|^{\zeta_t})\rho_t)|\mathcal{F}_t) \leq U'(0)E(|x| + (1+|x|^{\zeta_t})\rho_t|\mathcal{F}_t),$$

which shows (40), as well as (39), and the arguments of Section 4 provide an optimal $\tilde{\xi}_{t+1}(x)$ which satisfies (42). Skipping through arguments of Section 6 (Proposition 6.4), it becomes clear that for establishing the differentiability of $U_t$, one needs the integrability of something like

$$\sup_{x,y \in [a,b]} U'_{t+1}(x + y\Delta S_{t+1})|\Delta S_{t+1}|.$$

(We have switched to dimension 1 and $x, y \in [a, b]$ without loss of generality.) The above quantity is smaller than

$$|\Delta S_{t+1}|U'(-p(|\Delta S_{t+1}|\mathcal{X})) \leq |1 + |\Delta S_{t+1}| + \mathcal{X}|^r,$$

where $p(x) \geq 0$ is a polynomial of $x$; $r > 0$ and $\mathcal{X} \in \mathcal{M}$ [see (44), (42) and (38)], and this is indeed integrable. So the arguments of Section 6 apply, (43) and (44) follow. Finally, $EU_0(c) < \infty$ is deduced from (39); the existence of $\phi^*(c)$ and $Q$ follow just like in Sections 5 and 6. $\square$

NOTE ADDED IN PROOF. With a different argument it is posible to get rid of the left-hand side inequality in (38).

REMARK 7.2. Analogous arguments show that if $S$ is bounded, (NA) holds with $\kappa_t = \kappa, \beta_t = \beta$ deterministic constants in Proposition 3.3 and either one of Assumptions 2.3, 2.5 holds, then Assumption 2.1 is true and there exists a *bounded* optimal strategy $\phi^*$.

We now demonstrate that, even under (NA), the expected utility maximization problem does not necessarily admit an optimal portfolio.

EXAMPLE 7.3. Define a strictly increasing concave function $U$ by setting

$$\begin{aligned} U(0) &= 0, \\ U'(x) &:= 1 + 1/n^2, \quad x \in (n-1, n], n \geq 1, \\ U'(x) &:= 3 - 1/n^2, \quad x \in (n, n+1], n \leq -1. \end{aligned}$$

Take

$$S_0 := 0, \quad P(S_1 = 1) = 3/4, \quad P(S_1 = -1) = 1/4.$$

One can calculate the expected utility of the strategy $\phi_1 := n, n \in \mathbb{Z}$ with initial capital $c = 0$:

$$\begin{aligned} EU(nS_1) &= \frac{3U(n) + U(-n)}{4} \\ &= \frac{1}{4}\left(3n + 3\sum_{j=1}^{n} 1/j^2 - 3n + \sum_{j=1}^{n} 1/j^2\right) = \sum_{j=1}^{n} 1/j^2, \quad n \geq 0; \end{aligned}$$

and

$$EU(nS_1) = \tfrac{1}{4}\left(9n + 3\sum_{j=1}^{-n} 1/j^2 - n + \sum_{j=1}^{-n} 1/j^2\right) = \sum_{j=1}^{-n} 1/j^2 + 2n, \qquad n < 0.$$

This utility tends to

$$\sum_{i=1}^{\infty} 1/i^2 = \pi^2/6$$

in an increasing way as $n \to \infty$. In fact, it is easy to see that the function

$$\phi_1 \to EU(\phi_1 S_1), \qquad \phi_1 \in \mathbb{R}$$

is increasing in $\phi_1$, so we may conclude that the supremum of the expected utilities is $\pi^2/6$, but it is not attained by any strategy $\phi_1$. It is clear that one can construct a similar example with $U'$ continuous and $U$ strictly concave. It would be interesting to find the minimal conditions on $U$ which assure the existence of an optimal portfolio under (NA) and Assumption 2.1.

REMARK 7.4. In Section 6 we have proven for bounded $S$ a certain "individual" version of the fundamental theorem of asset pricing: absence of arbitrage implies that an agent of subjective utility $U$ finds an equivalent martingale measure computed from his or her optimal investment strategy. In the light of Proposition 7.1, we might relax the assumption on $S$. An interesting special case is when $U(x) := x, x \leq 0$ and otherwise $U$ and $S$ satisfy the assumptions of Proposition 7.1. Then we get a martingale measure with $dQ/dP$ bounded [due to the fact that $U'(x) \leq 1$]. That is, in this particular model class we have reproved the result of Dalang, Morton and Willinger [4].

REMARK 7.5. Using Proposition 7.1 with a suitable $U$ such that $U(x) = x$, $x \geq 0$, one obtains a martingale measure $Q \sim P$ such that $dQ/dP \geq h$ for some constant $h > 0$. This result seems to be new and did not follow from the generic functional analytic approach to the construction of martingale measures (i.e., separation theorems). Note that, just like in Remark 7.4, we rely on the fact that a larger class of utility functions is allowed, see Remark 2.9.

REMARK 7.6. We finally note that an optimal strategy exists under conic portfolio constraints too. Let $C$ be a fixed polyhedral cone in $\mathbb{R}^d$. If we admit only strategies satisfying $\phi_t \in C$ for all $t$ and define (NA), $\Phi$, and so on in the respective manner, Theorem 2.7 remains true. Modification is required only in the Fatou lemma arguments of Lemmas 4.8 and 4.9. The choice $C := \mathbb{R}_+^d$ corresponds to forbidden short sales. In this case, the argument of Theorem 6.2 provides a measure $Q \sim P$ such that $S$ is a $Q$-supermartingale. The fact that (NA) under short sales constraints is equivalent to the existence of an equivalent supermartingale-measure was first noticed in [16], see also [2, 22, 24].

## APPENDIX

Let $\mathcal{H} \subset \mathcal{F}$ be a $\sigma$-algebra containing $P$-zero sets. An $\mathcal{H}$-measurable *random set* $D$ is an element of $\mathcal{H} \otimes \mathcal{B}(\mathbb{R}^d)$, where $\mathcal{B}(\mathbb{R}^d)$ denotes the Borel sets of $\mathbb{R}^d$. A *random affine subspace* $D$ is $\mathcal{H}$-measurable random set such that $D(\omega)$ is an affine subspace of $\mathbb{R}^d$ for each $\omega$.

Let $Y$ be a $d$-dimensional random variable and $\mu(\cdot, \omega) := P(Y \in \cdot | \mathcal{H})$ a regular version of its conditional distribution. Let $D(\omega)$ be the smallest affine subspace of $\mathbb{R}^d$ containing the support of $\mu(\cdot, \omega)$.

PROPOSITION A.1. *$D$ is an $\mathcal{H}$-measurable random affine subspace.*

PROOF. This is only a sketch. We begin by showing that $\operatorname{supp} \mu(\cdot, \omega)$ or, equivalently, its complement $\operatorname{supp}^C \mu(\cdot, \omega)$ is a random set. Let $\mathcal{G}$ be a countable base for the topology of $\mathbb{R}^d$. Then

$$\operatorname{supp}^C \mu(\cdot, \omega) := \bigcup \{G \in \mathcal{G} : \mu(G, \omega) = 0\},$$

which proves the assertion. Actually, $Z(\omega) := \operatorname{conv}(\operatorname{supp} \mu(\cdot, \omega))$ is a random set, where $\operatorname{conv}(\cdot)$ denotes closed convex hull. This is based on the existence of a sequence of random variables that is dense in the random set $\operatorname{supp} \mu(\cdot, \cdot)$, which is provided by Theorem III. 22 on page 74 of [3]. Now a simple argument (Theorem III.40 on page 87 of [3]) shows that the closed convex hull is, indeed, a random set.

Take a measurable selector $\nu(\omega)$ of $Z(\omega)$. Then the random set $Z - Z$ contains $0$ in its relative interior,

$$\left[ \bigcup_{n \in \mathbb{N}} \{nz : z \in Z - Z\} \right] + \nu(\omega)$$

clearly equals $D(\omega)$, which proves the proposition. $\square$

LEMMA A.2. *Let $\eta_n : \mathbb{R} \times \Omega \to \mathbb{R}^d$ be a sequence of $\mathcal{B}(\mathbb{R}) \otimes \mathcal{H}$-measurable functions such that for almost all $\omega$,*

$$\forall x \qquad \liminf_{n \to \infty} |\eta_n(x, \omega)| < \infty.$$

*Then there is a sequence $n_k$ of $\mathcal{B}(\mathbb{R}) \otimes \mathcal{H}$-measurable $\mathbb{N}$-valued functions, $n_k < n_{k+1}, k \in \mathbb{N}$ such that almost surely $\tilde{\eta}_k(x, \omega) := \eta_{n_k}(x, \omega)$ converges for all $x$ to some $\tilde{\eta}(x, \omega)$ as $k \to \infty$. To put it more concisely, there is a convergent random subsequence.*

PROOF. This is just a variant of Lemma 2 in [17]. □

LEMMA A.3. *Let $C(z,\omega), z \in [0,1]^2$ be continuously differentiable for almost all $\omega$ and measurable for any fixed $z$ such that $E\sup_{z\in[0,1]^2}|C(z,\omega)| < \infty$ and $E[\sup_z |\partial_1 C(z,\omega)| + \sup_z |\partial_2 C(z,\omega)|] < \infty$. Then there is a version $H(z,\omega)$ of*

$$z \to E(C(z,\omega)|\mathcal{H}),$$

*which is almost surely a continuously differentiable function.*

PROOF. Regard $C$ as a random element of the Banach space $\mathbf{C}^1([0,1]^2)$ of continuously differentiable functions equipped with the norm

$$\|f\| := \sup_{z\in[0,1]^2} |f(z)| + \sup_{z\in[0,1]^2} |\partial_1 f(z)| + \sup_{z\in[0,1]^2} |\partial_2 f(z)|,$$

and the corresponding Borel-field. It is easy to see that $C$ is measurable in this sense. Then the assertion follows from Proposition V-2-5 of [23]. □

**Acknowledgments.** The research of this paper was carried out during M. Rásonyi's visits in the framework of IMPAN-BC Centre of Excellence; he wishes to express his gratitude for the invitation and the hospitality. The authors thank I. V. Evstigneev and P. Jaworski for their valuable comments.

COMPUTER AND AUTOMATION INSTITUTE
HUNGARIAN ACADEMY OF SCIENCES
1111 BUDAPEST
KENDE UTCA 13-17
HUNGARY
E-MAIL: rasonyi@sztaki.hu

INSTITUTE OF MATHEMATICS
POLISH ACADEMY OF SCIENCES
ŚNIADECKICH 8
00-950 WARSAW
POLAND
E-MAIL: stettner@impan.gov.pl